\documentstyle{article}

\newcommand{\Frac}[2]{\displaystyle\frac{#1}{#2}}
\newcommand{\Sum}{\displaystyle\sum}
\newcommand{\F}{\mbox{${\rm I} \! {\rm F}$}}
\newcommand{\dem}{\par{\em Proof\/: }\\ \noindent }
\newcommand{\findemo}{$\;\;\Box$\\}
\newfont{\frack}{eufm10}

\newtheorem{defi}{Definition}
\newtheorem{lem}{Lemma}
\newtheorem{thm}{Theorem}
\newtheorem{nota}{Remark}
\newtheorem{algor}{Algorithm}
\newtheorem{ejplo}{Example}

\begin{document}

\title{Decoding Algebraic Geometry codes by a key equation}
\author{J. I. Farr\'{a}n\thanks{Partially supported by DIGICYT PB94-1111-C02-01}}
\date{November 21, 1998}
\maketitle

\begin{abstract}

A new effective decoding algorithm is presented for arbitrary 
algebraic-geometric codes on the basis of solving a generalized 
key equation with the majority coset scheme of Duursma. 
It is an improvement of Ehrhard's algorithm, since 
the method corrects up to the half of the Goppa distance with 
complexity order ${\cal O}(n^{2.81})$\/, and with no further 
assumption on the degree of the divisor $G$\/. 


\end{abstract}

\section{Introduction}

Decoding algebraic-geometric codes (AG codes in short) in an 
effective way can be done by means of solving a key equation, 
generalizing the ideas of the Berlekamp-Massey algorithm for BCH 
codes or the Euclidean algorithm for classical Goppa codes (see \cite{Blk}). 
In the original version of Porter, Shen and Pellikaan (see \cite{PSP}), 
only one-point codes with further assumptions on the curve were decoded, 
but the main ideas of the method can be extended for arbitrary 
curves and AG codes with Ehrhard's version of the key equation. 
Nevertheless, this algorithm does not correct up to the Goppa distance, 
but the complexity is only ${\cal O}(n^{3})$ (more details in \cite{EhrTh}). 
Our aim is to include in this method a majority scheme which generalizes 
the ideas of Feng and Rao for one point codes (see \cite{FR}), 
together with giving an improvement of the complexity 
by using the new methods given in \cite{Strassen} to solve linear equations. 
Thus, the algorithm that we propose improves 
both the decoding capacity and the complexity 
without losing the generality of its application to arbitrary AG codes. 
It uses the majority coset decoding scheme, which was introduced by 
Duursma, with the only further assumption that there is an extra 
rational point in the curve which is not used in the construction 
of the codes (more details in \cite{MCD}). This hypothesis is actually 
a weakening of the assumptions required by Porter's method. 

In section {\bf 2} we rewrite Ehrhard's key equation in a way 
that is closer to the original ideas of Porter, Shen and Pellikaan, 
in order to show the explicit connection between both works. 
Afterwards, we summarize in section {\bf 3} the main ideas of 
Duursma's majority coset scheme, in order to give in section {\bf 4} 
an algorithm which includes the above majority scheme in the key equation, 
so that one can increase the error capacity without the assumption 
$deg\,G\geq 6g-2\tau-2$, where $\tau$ is the gonality of the curve, 
which is required in Ehrhard's algorithm given in \cite{EhrAlg} 
(see also \cite{DuurTh} for further details). 
In the paper, we fix a non-singular absolutely irreducible 
projective algebraic curve $\chi$ defined over $\F_{q}$ 
and rational points $P_{1},\ldots,P_{n}$ of $\chi$\/.

\section{Key equation and decoding}

Let $G$ be a rational divisor whose 
support is disjoint to $D=P_{1}+\ldots+P_{n}$. Assume that 
$2g-2<deg\,G<n+g$, and consider the code $C=C_{\Omega}(D,G)$, 
that is the image of the linear injective map 
$$res_{D}\;:\;\Omega(G-D)\rightarrow\F_{q}^{n}$$
$$\eta\mapsto(res_{P_{1}}(\eta),\ldots,res_{P_{n}}(\eta))$$
with dimension $k\geq n-deg\,G+g-1$ and minimum distance $d\geq 
d^{\ast}=deg\,G+2-2g$, where $g$ is the genus of the curve. 
In the sequel, we fix a divisor $G^{\ast}$ with $\ell(G^{\ast})=0$ 
and $G\geq G^{\ast}\,$. In order to decode $C$, 
we will give a result for preparation.

\begin{lem}

There exists a vector space $V$ of differential forms such that 
$\Omega(G-D)\subseteq V$ and $res_{D}\;:\;V\rightarrow\F_{q}^{n}$ 
is an isomorphism. 

\end{lem}

\dem

Since $\Omega(G-D)\subseteq\Omega(G^{\ast}-D)$\/, it suffices to prove 
that $res_{D}$ is surjective on $\Omega(G^{\ast}-D)$, because it is 
injective on $\Omega(G-D)$. But the kernel of $res_{D}$ 
considered on $\Omega(G^{\ast}-D)$ is $\Omega(G^{\ast})$; hence the 
rank is $i(G^{\ast}-D)-i(G^{\ast})=
deg\,G^{\ast}-deg\,(G^{\ast}-D)=n$, because of the Riemann-Roch formula. 

\vspace{.2cm}
\findemo

\begin{nota}

In the sequel we fix an arbitrary differential form $\eta\neq 0$ and write 
$K=(\eta)$\/. Then for any rational divisor $H$ consider the isomorphism 
$${\cal L}(K-H)\rightarrow\Omega(H)$$
given by 
$$f\mapsto f\eta$$
This map is compatible with inclusions and restrictions, and so 
the inclusions $\Omega(G-D)\subseteq V\subseteq\Omega(G^{\ast}-D)$
give the corresponding ${\cal L}(K+D-G)\subseteq U\subseteq
{\cal L}(K+D-G^{\ast})$, where the map $f\mapsto res_{D}(f\eta)$ is 
an isomorphism from $U$ onto $\F_{q}^{n}$. Denote the inverse of this 
last map by ${\bf y}\mapsto h_{\bf y}$, i.e. $h_{\bf y}$ is the unique 
element in $U$ such that $res_{D}(h_{\bf y}\eta)={\bf y}$. 

\end{nota}

Because of the bijection $C\stackrel{\sim}{\leftrightarrow}
{\cal L}(K+D-G)$ given by ${\bf y}\leftrightarrow h_{\bf y}\,$, 
the {\em decoding problem} can be obviously described as follows: 
\begin{quote}
$(\ast)\;\;$
\begin{em} Given ${\bf y}\in\F_{q}^{n}\,$, find a function 
$h_{\bf c}\in{\cal L}(K+D-G)$ such that $h_{\bf e}\eta$ 
has a minimal number of poles in $sup\,(D)$\/, 
where $h_{\bf e}=h_{\bf y}-h_{\bf c}\,$\end{em}. 
\end{quote}
This problem will be solved by the following definition and results.

\begin{defi}

Given an arbitrary divisor $F$, a solution of the key 
equation for the received word ${\bf y}$ (related to $F$) 
is a triple $(f,q,r)\in({\cal L}(F)\setminus\{0\})
\times{\cal L}(K+F+D-G)\times{\cal L}(K+F-G^{\ast})$ 
such that $f\,h_{\bf y}=q+r\,$. 

\end{defi}

Notice that this definition means that $h_{\bf y}=\Frac{q}{f}+\Frac{r}{f}$ 
and $h_{\bf c}=\Frac{q}{f}\in{\cal L}(K+D-G)$\/. Thus, what we need to solve 
the decoding problem is giving conditions so that $h_{\bf e}=\Frac{r}{f}$ 
has few poles in $sup\,(D)$\/. This is done by the following theorem.

\begin{thm}[Decoding theorem] Let ${\bf y}={\bf c}+{\bf e}$, where 
${\bf c}\in C$\/. Then: 

\begin{description}

\item[1. ] If ${\cal L}(F-D_{\bf e})\neq 0$, then there exists a solution 
of the key equation. 

\item[2. ] If $deg\,F+wt({\bf e})<d^{\ast}$, then any solution 
$(f,q,r)$ of the key equation satisfies 
$$res_{D}\left(\frac{q\eta}{f}\right)={\bf c}\;\;\;\;{\rm and}\;\;\;\;
res_{D}\left(\frac{r\eta}{f}\right)={\bf e}$$

\end{description}

\end{thm}

\dem

\begin{enumerate}

\item Take a non-zero function 
$f\in{\cal L}(F-D_{\bf e})\subseteq{\cal L}(F)$. Then 
$(f\,h_{\bf c})\geq G-D-F-K$, $(f\,h_{\bf e})\geq -F+D_{\bf e}+
G^{\ast}-D_{\bf e}-K
=G^{\ast}-F-K$ and $f\,h_{\bf y}=f\,h_{\bf c}+f\,h_{\bf e}$; hence 
the triple $(f,f\,h_{\bf c},f\,h_{\bf e})$ is a solution of the 
key equation. 

\item Denote by $D_{\bf e}$ the divisor of poles of $h_{\bf e}\eta$ 
in the support of $D$\/. Let $(f,q,r)$ be a solution of the key equation 
and set $\varphi\doteq r-f\,h_{\bf e}=f\,h_{\bf c}-q$\/. 
One can estimate the following divisors: 
$$K+(r-f\,h_{\bf e})\geq\min\{G^{\ast}-F,G^{\ast}-F-D_{\bf e}\}=
G^{\ast}-F-D_{\bf e}$$
and 
$$K+(f\,h_{\bf c}-q)\geq G-F-D$$
what means that $\varphi\in{\cal L}(K+F+D_{\bf e}-G^{\ast})\cap
{\cal L}(K+F+D-G)={\cal L}(K+F+D_{\bf e}-G)=0$, since by assumption 
$deg\,(K+F+D_{\bf e}-G)=2g-2+deg\,(F)+wt({\bf e})-deg\,(G)<0$. 
Hence $\varphi=r-f\,h_{\bf e}=f\,h_{\bf c}-q=0$, what yields the theorem. 

\end{enumerate}

\vspace{.2cm}
\findemo

Assume from now on that ${\cal L}(F-D_{\bf e})\neq 0$ and 
$deg\,F+wt({\bf e})<d^{\ast}$ (notice that both assumptions are 
satisfied if $wt({\bf e})\leq\nu$ and $deg\,(F)=\nu+g$\/, where 
$\nu\doteq\left\lfloor\Frac{d^{\ast}-g-1}{2}\right\rfloor$\/; 
that is, when there are few errors and $F$ is small). 
Thus, for a fixed ${\bf y}\in\F_{q}^{n}$ define the linear map 
$$\begin{array}{ccccc}
\varepsilon_{\bf y}&:&{\cal L}(F)&\rightarrow&{\cal L}(K+F+D-G^{\ast})\\
\,&\,&f&\mapsto&f\,h_{\bf y}\end{array}$$
Since $deg\,(G-F)>deg\,G-d^{\ast}=2g-2$, one has 
${\cal L}(K+F+D-G)\cap{\cal L}(K+F-G^{\ast})={\cal L}(K+F-G)=0$, and hence 
there exists a vector space $W$ such that 
$${\cal L}(K+F+D-G^{\ast})={\cal L}(K+F+D-G)\oplus{\cal L}(K+F-G^{\ast})
\oplus W$$

Denoting by $\pi_{W}$ and $\pi^{\ast}$ the natural 
projections onto $W$ and ${\cal L}(K+F-G^{\ast})$ respectively, notice 
that the key equation means that $\varepsilon_{\bf y}(f)$ has a 
null projection onto $W$. Therefore, if there exists a codeword 
${\bf c}$ satisfying $wt({\bf y}-{\bf c})\leq t$, where 
$0<t\leq\left\lfloor\Frac{d^{\ast}-g-1}{2}\right\rfloor$ is fixed, 
one can compute the error vector with the following algorithm, 
where a suitable basis for every above function space is assumed 
to be previously calculated. Such bases can be computed 
by means of Brill-Noether algorithm (see \cite{HachLB}).

\begin{algor}[${\cal K}_{G}(F)$]

\begin{description}

\item[] 

\item[1. ] Compute a matrix for the linear map $\varepsilon_{\bf y}$\/. 

\item[2. ] Find a non-zero function 
$f\in ker\,(\pi_{W}\circ\varepsilon_{\bf y})$\/. 

\item[3. ] Compute $r=\pi^{\ast}(\varepsilon_{\bf y}(f))$\/. 

\item[4. ] Compute ${\bf e}=res_{D}\Frac{r\eta}{f}$, checking that 
${\bf y}-{\bf e}\in C$ and $wt({\bf e})\leq\left\lfloor
\Frac{d^{\ast}-1}{2}\right\rfloor$\/. 

\end{description}

\end{algor}

Notice that most of the calculations in this algorithm are concentrated 
in the first two steps, and thus its complexity is that of solving 
linear equations (see \cite{EhrTh}). Also notice that the algorithm may fail 
in the second or forth steps if the number of comitted errors is greater 
than the bound $\left\lfloor\Frac{d^{\ast}-g-1}{2}\right\rfloor$\/, 
and hence it cannot correct in general up to the half of the Goppa 
distance. In order to do it, we can use a majority voting scheme, 
what will be explained in the next sections.

\begin{nota}

We show now how the above results generalize those of 
Porter, Shen and Pellikaan, and why they are stronger. 
Following the notations from \cite{PSP}, the original algorithm works 
with the codes $C=C_{\Omega}(D,E-\mu P)$, where $E$ is the divisor of 
zeros of a function $h\in K_{\infty}(P)$ without zeros in $sup\,(D)$\/, 
$K_{\infty}(P)$ being the ring of those functions having poles only 
at $P$\/, $P$ being a rational point distinct from $sup\,(D)$\/, 
and where $\mu$ is a positive integer. 
In this case, we can obviously take $G^{\ast}=-\mu P$\/. 
For the sake of simplicity, assume that there exists 
a differential form $\eta$ such that $(\eta)=(2g-2)P$\/. 

Firstly, from the isomorphism given by lemma {\bf 1} we obtain 
a basis $\varepsilon_{1},\ldots,\varepsilon_{n}$ of $V$ such that 
$res_{D}(\varepsilon_{1}),\ldots,res_{D}(\varepsilon_{n})$ is the 
canonical basis of $\F_{q}^{n}\,$. 
Then, Porter defines a \lq\lq syndrome function" by 
$$S_{\bf y}\cdot\eta=
\Sum_{j=1}^{n}y_{j}\left(1-\Frac{h}{h(P_{j})}\right)\varepsilon_{j}$$
Notice that $S_{\bf y}\in K_{\infty}(P)$\/, 
$S_{\bf y}\equiv h_{\bf y}\;\;(mod\;h)$ and 
$-\upsilon_{P}(S_{\bf y})\leq m+2g-1$, where $m\doteq-\upsilon_{P}(h)$\/. 
On the other hand, Porter's result to decode $C$ can be rewritten 
as follows (see \cite{EhrTh} for further details): 
\begin{quote}
If there is an integer $t$ such that $t+wt({\bf e})<d^{\ast}\,$ 
and functions $f,q,r\in K_{\infty}(P)$ satisfying 
$-\upsilon_{P}(r)\leq t+2g-2+\mu$, $-\upsilon_{P}(f)\leq t$ 
and the \lq\lq polynomial key equation" 
$$f\,S_{\bf y}=gh+r$$
then $h_{\bf e}=\Frac{r}{f}$\/. 
\end{quote}

Such triples $(f,g,r)$ are called \lq\lq valid solutions" 
in \cite{PSP}. Thus, by taking $K=(\eta)=(2g-2)P$ and $F=tP$\/, 
one has $f\in{\cal L}(F)$ and $r\in{\cal L}(K+F-G^{\ast})$\/, 
and hence this is a particular case of our method
\footnote{In particular, the condition of being minimal for a 
valid solution can be dismissed from the results of Porter. }. 
Moreover, one obtains ${\bf e}=res_{D}\Frac{r\eta}{f}$ where $f$ 
has few zeros for $F$ \lq\lq small" (because of $(f)+F\geq 0$) 
and thus, for a suitable choice of $t$ and $wt({\bf e})$\/, 
$\,\Frac{r\eta}{f}$ has a minimal number of poles in $sup\,(D)$\/, 
according to the formulation ($\ast$) of the decoding problem. 
This is actually the underlying idea of Porter, which was 
carried out by a \lq\lq row reduction process" at a certain 
resultant matrix, but of course it can also be done by simple 
techniques of linear algebra, as we have explained above. 

Thus, the results of our paper are stronger than the originals, 
since they work with an arbitrary divisor $G$ and we do not 
require any special differential form $\eta$ or rational 
function $h$\/, what is actually a very strong restriction. 
Moreover, one obtains a quite similar formula to compute the error 
just from $h_{\bf y}$\/, without the need of the syndrome $S_{\bf y}$\/. 

\end{nota}

\section{Majority coset decoding}

This section is abstracted from \cite{MCD}. 
Assume that there exists a rational point $P_{\infty}\notin sup\,(D)$\/, and 
let $H_{1}$ be a rational divisor whose support is disjoint to $sup\,(D)$\/. 
Set $H_{0}=H_{1}-P_{\infty}$ and $H_{2}=H_{1}+P_{\infty}\,$. For $i=0,1,2\;$, 
let $C_{i}=C_{\Omega}(D,H_{i})$ and $d_{i}^{\ast}=deg\,(H_{i})+2-2g$\/. 
One obviously has $C_{0}\supseteq C_{1}\supseteq C_{2}\,$. 

For an error vector ${\bf e}$ such that $wt({\bf e})\leq(d_{1}^{\ast}-1)/2$ 
we want to solve the following problem: 
\begin{quote}\begin{em} 
Given ${\bf y}_{1}$ with ${\bf y}_{1}-{\bf e}\in C_{1}$, 
finding ${\bf y}_{2}$ such that ${\bf y}_{2}-{\bf e}\in C_{2}\,$. 
\end{em}\end{quote} 
Such a problem is called \begin{em} coset decoding procedure 
related to the extension $C_{1}\supseteq C_{2}\,$,\/ \end{em} 
where we obviously can assume that $C_{1}\neq C_{2}\,$. 

Thus, for a given ${\bf y}\in\F_{q}^{n}$ 
and for any rational function $h$ without poles in $sup\,(D)$\/, 
one defines the syndrome $S_{\bf y}(h)$ by the expression 
$$S_{\bf y}(h)\doteq\sum_{j=1}^{n}y_{j}h(P_{j})\in\F_{q}$$
which is linear with respect to both ${\bf y}$ and $h$\/. 

It is very easy to prove that the syndrome is a {\em coset invariant}\/, 
i.e. $S_{\bf y}(h)=S_{\bf e}(h)$ for all $h\in{\cal L}(H_{i})$ 
if and only if ${\bf y}-{\bf e}\in C_{i}$, for $i=0,1,2\,$. 
Hence, ${\bf y}\in C_{i}$ if and only if 
$S_{\bf y}(h)=0$ for all $h\in{\cal L}(H_{i})$\/. 

On the other hand, for an arbitrary divisor $F$ defined over $\F_{q}$ 
and $i=0,1,2\,$, one defines the kernels $K_{i}(F)$ 
associated to the error vector ${\bf e}$ by 
$$K_{i}(F)\doteq\{f\in{\cal L}(F)\;|\;S_{\bf e}(f\cdot g)=0,\;\;\forall 
g\in{\cal L}(H_{i}-F)\}$$
All the vector spaces 
$K_{1}(F+P_{\infty})/K_{0}(F)$\/, $K_{0}(F)/K_{1}(F)$\/, 
${\cal L}(H_{1}-F)/{\cal L}(H_{1}-F-P_{\infty})$\/, 
$K_{1}(F+P_{\infty})/K_{2}(F+P_{\infty})$ and 
$K_{2}(F+P_{\infty})/K_{1}(F)$ have dimension at most one. 
Thus, we are interested in the following conditions: 
$$\begin{array}{lcl}
{\bf (A1)}\;K_{1}(F+P_{\infty})\neq K_{0}(F)&\;\;&
{\bf (B1)}\;K_{1}(F+P_{\infty})=K_{2}(F+P_{\infty})\\
{\bf (A2)}\;K_{0}(F)=K_{1}(F)&\;\;&
{\bf (B2)}\;K_{2}(F+P_{\infty})\neq K_{1}(F)\\
{\bf (A3)}\;{\cal L}(H_{1}-F)\neq{\cal L}(H_{1}-F-P_{\infty})&\;\;&\,
\end{array}$$

Define the conditions $(A)\Leftrightarrow(A1)\wedge(A2)\wedge(A3)$ 
and $(B)\Leftrightarrow(B1)\wedge(B2)$\/. Since one has 
$(A1)\wedge(B1)\Leftrightarrow(A2)\wedge(B2)$\/, the conditions 
$(A)$ and $(B)$ are equivalent to $(A1)$, $(A3)$ and $(B1)$\/. 

It follows from \cite{MCD} (sections II and III) that if $(A)$ and $(B)$ 
are satisfied, then the {\em coset decoding procedure}\/ can be 
implemented by the following algorithm, 
where $D$ and $P_{\infty}$ are fixed.

\begin{algor}[${\cal C}_{H_{1}}(F)$]

\begin{description}

\item[] 

\item[] $Input:={\bf y}_{1}\,$. 

\item[] If $C_{1}=C_{2}$ then ${\bf y}_{2}={\bf y}_{1}\,$ else: 

\begin{itemize}

\item Find ${\bf c}_{0}\in C_{1}\setminus C_{2}\,$. 

\item Find $f\in K_{1}(F+P_{\infty})\setminus K_{0}(F)$\/. 

\item Find $g\in{\cal L}(H_{1}-F)\setminus{\cal L}(H_{1}-F-P_{\infty})$\/. 

\item Compute $\lambda=S_{{\bf y}_{1}}(fg)/S_{{\bf c}_{0}}(fg)$\/. 

\item Set ${\bf y}_{2}={\bf y}_{1}-\lambda{\bf c}_{0}\,$. 

\end{itemize}

\item[] $Output:={\bf y}_{2}\,$. 

\end{description}

\end{algor}

Unfortunately we are not able in practice to check the condition $(B)$\/, 
since $K_{2}(F+P_{\infty})$ is not known from the received word ${\bf y}$. 
This problem can be solved by means of a {\em majority voting}\/, 
on the basis of the following result due to Duursma 
(see \cite{MCD} for further details). 

\vspace{.3cm}

\begin{thm}[Main theorem]

Let $C_{0}\supseteq C_{1}\supseteq C_{2}$ be the extension of codes 
given by $C_{i}\doteq C_{\Omega}(D,H_{i})$\/, where $H_{1}$ 
has disjoint support with $D$\/, $H_{0}\doteq H_{1}-P_{\infty}$ and 
$H_{2}\doteq H_{1}+P_{\infty}\,$. Assume that the genus is $g\geq 1$, 
and take numbers $t,r\geq 0$ such that 
$2t+r+1\leq d_{1}^{\ast}\doteq deg\,H_{1}+2-2g$\/. 
Take an arbitrary divisor $F_{0}$ of degree $t$, 
and define $F_{i}\doteq F_{0}+iP_{\infty}$ for $i=1,\ldots,2g-1$. For an 
error vector ${\bf e}$ with weight $wt({\bf e})\leq t$, define: 
$$I\doteq\{r,r+1,\ldots,2g-2\}$$
$$T\doteq\{i\in I\;|\;(A)\wedge(B)\;\;{\rm hold}\;\;{\rm for}\;\;F=F_{i}\}$$
$$F\doteq\{i\in I\;|\;(A)\wedge\neg(B)\;\;{\rm hold}\;\;{\rm for}\;
\;F=F_{i}\}$$
Then at least one of the following conditions holds: 
$$\begin{array}{ccc}
{\bf (i)}&\hspace{1.cm}&
{\cal L}(H_{1}-F_{2g-1}-D_{\bf e}-rP_{\infty})\neq 0\\
{\bf (ii)}&\hspace{1.cm}&
{\cal L}(F_{r}-D_{\bf e})\neq 0\\
{\bf (iii)}&\hspace{1.cm}&\sharp T>\sharp F
\end{array}$$

\end{thm}

In the last section we will see how to apply this majority scheme 
in order to improve the correction capacity of the decoding algorithm 
by solving the Ehrhard's key equation up to the half of the Goppa distance. 
The so obtained procedure is thus the best possible one by solving a key 
equation, looking at the generality and the capacity of the algorithm.

\section{Decoding by a key equation with majority voting}

Let $C=C_{\Omega}(D,G)$ be a {\em strongly algebraic-geometric code}\/, 
i.e. such that $2g-2<deg\,(G)<n$\/. For our purpose, we can assume 
that $g>0$, since otherwise the key equation corrects $C$ up to the 
half of the Goppa distance and we do not need any majority voting. 

Consider successive divisors $G_{r}=G+rP_{\infty}$, for $r=0,1,\ldots,g$\/. 
Notice that for any such divisor $G_{r}$ one has $2g-2<deg\,(G_{r})<n+g$\/, 
and thus all these divisors are in the situation of the first paragraph 
in section {\bf 2}. On the other hand, 
take $t\doteq\left\lfloor\Frac{d^{\ast}-1}{2}\right\rfloor$\/, 
where $d^{\ast}\doteq deg\,(G)+2-2g$\/, and assume $t>0$. 
Take then a divisor $F_{0}$ with degree $t$ and 
set $F_{i}\doteq F_{0}+iP_{\infty}$ for $i=1,\ldots,2g-1\,$. 

Thus we can consider the following algorithm, 
which brings together the methods of Ehrhard and Duursma. 
In the algorithm, the main idea is that the conditions {\bf (i)} and 
{\bf (ii)} given by theorem {\bf 2} allows us to get the error vector 
by means of a key equation for some suitable $G$ and $F$, and otherwise 
the condition {\bf (iii)} provides us with a majority test to solve 
the coset decoding problem and decrease the size of the code. 
We assume that bases for the involved function and differential spaces are 
previously calculated together with the spaces $U,V,W$ as in section {\bf 2}, 
for all of the possible cases when algorithm {\bf 1} is applied.

\begin{algor}[${\cal D}_{G}(F_{0})$]

\begin{description}

\item[] 

\item[] $Input:={\bf y}\in\F_{q}^{n}\,$. 

\item[] Set ${\bf y}_{1}={\bf y}\,$. 

\item[] From $r=0$ to $r=g$ do: 

\begin{itemize}

\item Set $H_{1}=G+rP_{\infty}\,$. 

\item If ${\cal K}_{G}(G-F_{2g-1})$ gets the error vector 
from ${\bf y}_{1}\,$, then return ${\bf e}$ and STOP. 

\item Otherwise, if ${\cal K}_{H_{1}}(F_{r})$ gets the error vector 
from ${\bf y}_{1}\,$, then return ${\bf e}$ and STOP. 

\item Otherwise, compute $I_{A}\doteq\{i=r,r+1,\ldots,2g-2\;|\;\,(A)\;\,
{\rm holds}\;\,{\rm for}\;\,F=F_{i}\}$, apply the coset decoding procedure 
${\cal C}_{H_{1}}(F_{i})$ for $i\in I_{A}$ with input ${\bf y}_{1}$ 
and get a vector ${\bf y}_{2}$ whose coset with respect to 
$C_{2}\doteq C_{\Omega}(D,H_{1}+P_{\infty})$ occurs most of the times. 

\vspace{.2cm}

\noindent Set ${\bf y}_{1}={\bf y}_{2}\,$ and NEXT $r$\/. 

\end{itemize}

\end{description}

\end{algor}

Notice that algorithm {\bf 1} is always applied to one of the 
divisors $G_{r}\,$. Thus, if we take a divisor $G^{\ast}$ such that 
$\ell(G^{\ast})=0$ and $G^{\ast}\leq G\leq G_{r}\,$, we can use the 
same divisor $G^{\ast}$ for all the involved key equations. 

Finally, since every functional code can be expressed as a differential 
code and vice versa, we can prove the following new result, 
which incooperates the Duursma's version of the majority voting scheme 
into the Ehrhard's version of the key equation.

\begin{thm}

Let $\chi$ be a non-singular absolutely irreducible projective 
algebraic curve defined over the finite field $\F_{q}$ with at least 
$n+1$ rational points. Let $C=C_{\Omega}(D,G)$ be an algebraic-geometric 
code with length $n$ such that $2g-2<deg\,(G)<n$\/. Let $F_{0}$ be any 
divisor with degree $t\doteq\left\lfloor\Frac{d^{\ast}-1}{2}\right\rfloor$\/, 
where $d^{\ast}\doteq deg\,(G)+2-2g$ is the Goppa distance of $C$\/. 
Then the algorithm ${\cal D}_{G}(F_{0})$ decodes $C$ up to $t$ errors 
with complexity ${\cal O}(n^{2.81})$\/. 

\end{thm}

\dem 

First of all, the condition $2t+r+1\leq d_{1}^{\ast}=deg\,(H_{1})+2-2g$ 
is satisfied by every divisor $H_{1}=G_{r}$ from $r=0$ to $r=g$\/, and 
for $t\doteq\left\lfloor(d^{\ast}-1)/2\right\rfloor$\/; 
thus we can apply theorem {\bf 2} in every step of the algorithm, 
provided $wt({\bf e})\leq t$\/. 

For a fixed $H_{1}=G_{r}\,$, if the condition {\bf (i)} 
${\cal L}(H_{1}-F_{2g-1}-D_{\bf e}-rP_{\infty})=
{\cal L}(G-F_{2g-1}-D_{\bf e})\neq 0$ holds 
together with $wt({\bf e})\leq t$\/, then the key equation 
${\cal K}_{G}(F)$ obtains the error vector for $F=G-F_{2g-1}\,$, 
since $deg\,F+wt({\bf e})<d^{\ast}$ and ${\cal L}(F-D_{\bf e})\neq 0$, 
and theorem {\bf 1} can be applied. 

In the same way, 
if the condition {\bf (ii)} ${\cal L}(F_{r}-D_{\bf e})\neq 0$ holds 
together with $wt({\bf e})\leq t$\/, then the key equation 
${\cal K}_{G_{r}}(F_{r})$ obtains the error vector, since 
$deg\,F_{r}+wt({\bf e})<deg\,(G_{r})+2-2g$ and 
${\cal L}(F_{r}-D_{\bf e})\neq 0$, and theorem {\bf 1} can also be applied. 

Otherwise, the condition {\bf (iii)} implies that the algorithm 
${\cal C}_{G_{r}}(F_{i})$ is correct for most of the \lq\lq candidates" 
$i\in I_{A}\,$, and we can carry on with the next step. 
Finally, for $r=g$ the condition ${\cal L}(F_{r}-D_{\bf e})\neq 0$ 
is always true and the algorithm stops at most in $g+1$ steps, 
if not too many errors occur. 

Notice that still the complexity of this algorithm is equivalent to solve 
a linear system of size $n$\/, since most of the computations come from 
either applications of the algorithm ${\cal K}_{G}(F)$ or finding a function 
in $K_{1}(F+P_{\infty})\setminus K_{0}(F)$ (more details in \cite{MCD}). 
Thus, the complexity is actually ${\cal O}(n^{2.81})$ 
\footnote{Nowadays there are even some improvements of this complexity. }, 
since solving linear equations can be done faster than Gaussian elimination 
(see for instance \cite{Strassen}). 

\vspace{.2cm}
\findemo

\begin{nota}

Notice that the complexity ${\cal O}(n^{2.81})$ is even better than 
the complexity of Sakata's algorithm ${\cal O}(n^{3-\frac{2}{r+1}})$ 
if the curve $\chi$ is embedded in an affine $r$\/-space with $r>10$ 
(what happens in the constructions of asymptotically good codes 
given in \cite{GStich}). Thus, general decoding methods which are 
based on solving linear equations are not so far from \lq\lq fast 
decoding" as they are supposed to (see \cite{HP} for a survey on decoding). 

\end{nota}

\begin{ejplo}

Consider the Klein quartic 
${\rm X}^{3}{\rm Y}+{\rm Y}^{3}{\rm Z}+{\rm Z}^{3}{\rm X}=0$ 
over $\F_{8}$. This curve has genus $g=3$ and $24$ rational points, 
namely, $Q_{0}=(1:0:0)$, $Q_{1}=(0:1:0)$ and $Q_{2}=(0:0:1)$ on the 
coordinate lines, and all the others are in the affine plane, 
namely, $P_{1},\ldots,P_{21}$ (see \cite{Hansen} for details). 
Set $H_{1}=G=4(Q_{0}+Q_{1}+Q_{2})$, $D=P_{1}+\ldots+P_{21}$ 
and define the code $C_{1}=C=C_{\Omega}(D,G)$, with parameters 
$[21,11,\geq 8]$. Consider the vector 
${\bf y}_{1}=(1,0,1,\alpha,0,\ldots,0))$ 
as a received word, where $\alpha\in\F_{8}$ 
satisfies $\alpha^{3}+\alpha+1=0$, 
and take the divisor $F_{0}=3P_{\infty}\,$, where $P_{\infty}=Q_{2}\,$. 
Notice that the correction capacity of our algorithm is $t=3$, 
whereas the key equation only corrects two errors. 

Thus, in the step $r=0$ one easily checks that the conditions {\bf (i)} 
and {\bf (ii)} from theorem {\bf 2} are not satisfied, and hence the key 
equation cannot correct this error. Then, one computes the set $I_{A}=\{3\}$ 
and applies ${\cal C}_{G}(F)$ to the only candidate $F=F_{3}\,$: 

\begin{itemize}

\item Take ${\bf c}=(\alpha,\alpha^{5},\alpha^{3},0,\alpha^{4},
\alpha^{2},\alpha^{6},1,0,1,1,0,\ldots,0)\in C_{1}\setminus C_{2}\,$. 

\item Take $f=\alpha^{3}+\Frac{{\rm Z}^{3}}{{\rm X}^{2}{\rm Y}}
\in K_{1}(F_{3}+P_{\infty})\setminus K_{0}(F_{3})$\/. 

\item Take $g=\Frac{\rm X}{\rm Y}
\in{\cal L}(G-F_{3})\setminus{\cal L}(G-F_{3}-P_{\infty})$. 

\item Compute $\lambda=S_{{\bf y}_{1}}(fg)/S_{{\bf c}}(fg)=\alpha^{3}$ 

\item Return ${\bf y}_{2}={\bf y}_{1}-\lambda{\bf c}=
(\alpha^{5},\alpha,\alpha^{2},\alpha,1,\alpha^{5},\alpha^{2},\alpha^{3},0,
\alpha^{3},\alpha^{3},0,\ldots,0)$\/. 

\end{itemize}

In this case we have no voting since there is an only candidate, and 
the above solution is the new ${\bf y}_{1}$ for the next step of the 
algorithm, which works in a smaller code, and go on until the key equation 
gets the error vector. 

\end{ejplo}

\begin{ejplo}

Consider now the Hermite curve 
${\rm Y}^{4}{\rm Z}+{\rm Y}{\rm Z}^{4}+{\rm X}^{5}=0$ 
over $\F_{16}\,$. It has $64$ affine rational points and only one point 
$P_{\infty}$ at infinity. Let $D=P_{1}+\ldots+P_{64}$, $G_{1}=23P_{\infty}$ and 
define the code $C=C_{\Omega}(D,G_{1})$, which is of type $[64,46,\geq 13]$\/. 
Consider then ${\bf y}_{1}=(\alpha^{12},\alpha^{4},\alpha^{7},\alpha^{8},
\alpha^{9},\alpha^{9},0,\ldots,0))$ as a received word, where $\alpha\in\F_{16}$ 
satisfies $\alpha^{4}+\alpha+1=0$, and take the divisor $F_{0}=6P_{\infty}\,$. 

Now for $r=0$ again {(i)} and {(ii)} do not hold, and one computes 
$I_{A}=\{1,2,3,5,7,8,9\}$\/. In this case, voting actually occurs and 
the procedure is equivalent to the algorithm of Feng and Rao (see \cite{MCD}). 

\end{ejplo}

\end{document}